\font\we=cmb10 at 14.4truept
\font\li=cmb10 at 12truept
\noindent
\centerline {\we Non-Abelian Class Field Theory for Riemann Surfaces}
\vskip 0.45cm
\centerline {\li Lin WENG}
\vskip 0.30cm
\centerline {\bf Graduate School of Mathematics, Nagoya University, Japan}
\vskip 0.50cm
\noindent
In this paper, using what we call a micro reciprocity law, we complete
Weil's program [W] for non-abelian class field theory of Riemann
surfaces.
\vskip 0.30cm
\noindent
{\li 1. Refined Structures for Tannakian Categories}
\vskip 0.30cm
Let {\bf T} be a Tannakian category with a fiber functor $\omega:{\bf
T}\to {\rm Ver}_{\bf C}$, where ${\rm Ver}_{\bf C}$ denotes the category
of finite dimensional {\bf C}-vector spaces. An object $t\in {\bf T}$
is called {\it reducible} if there exist non-zero objects $x,y\in {\bf T}$
such that $t=x\oplus y$. An object is called {\it irreducible} if it is
not reducible. If moreover every object $x$ of {\bf T} can be written
uniquely as a sum of irreducible objects $x=x_1\oplus
x_2\oplus\dots\oplus x_n$, then {\bf T} is called a {\it unique
factorization} Tannakian category. Usually, we call $x_i$'s the
{\it irreducible components} of $x$.

A Tannakian subcategory {\bf S} of a unique factorization Tannakian
category {\bf T} is called {\it completed} if for $x\in {\bf S}$,
all its irreducible components $x_i$'s in {\bf T} are also in {\bf S}.
{\bf S} is called {\it finitely generated} if as a Tannakian category, it
is generated by finitely many objects. Moreover,  {\bf S}  is called {\it
finitely completed} if (a) {\bf S} is finitely generated; (b) {\bf S} is
completed; and (c) ${\rm Aut}^{\otimes}\omega\big|_{\bf S}$ is a finite
group.
\vskip 0.30cm
\noindent
{\li 2. An Example}
\vskip 0.30cm
Let $M$ be a compact Riemann surface of genus $g$. Fix an effective
divisor $D=\sum_{i=1}^Ne_iP_i$ with $e_i\in {\bf Z}_{\geq 2}$ once for
all. For simplicity, in this note we always assume that $(M;D)\not=({\bf
P}^1;e_1P_1)$, or $({\bf P}^1;e_1P_1+e_2P_2)$ with $e_1\not=e_2$. (These
cases may be easily treated.)

By definition a parabolic semi-stable bundle

\centerline {$\Sigma:=(E=:E(\Sigma);P_1,\dots,P_N; Fil (E|{P_1}),\dots,
Fil(E|_{P_N});a_{11},\dots, a_{1r_1};\dots; a_{N1},\dots, a_{Nr_N})$}

\noindent
of parabolic degree 0 is called a {\it GA bundle} over $M$ along $D$ if
(i) the parabolic weights are all rational, i.e.,
$a_{ij}\in {\bf Q}\cap [0,1)$; (ii) there exist
$\alpha_{ij}\in {\bf Z}, \beta_{ij}\in {\bf Z}_{>0}$ such that (a)
$(\alpha_{ij},\beta_{ij})=1;$ (b) $a_{ij}=\alpha_{ij}/\beta_{ij};$ and (c)
$\beta_{ij}|e_i$, for all $i,j$. Denote by
$[\Sigma]$ the Seshadri equivalence class associated with $\Sigma$.
Moreover, for $[\Sigma]$, define $\omega_D([\Sigma])$ as $E({\rm
Gr}(\Sigma))|_P$, i.e., the fiber of the bundles associated with
Jordan-H\"older graded parabolic bundle of $\Sigma$ at a fixed $P\in
M^0:=M\backslash |D|=M\backslash\{P_1,\dots,P_N\}$.
\vskip 0.30cm
\noindent
{\bf Proposition.} {\it With  the same notation as above, put ${\cal
M}(M;D):=\{[\Sigma]:\Sigma\ {\rm is\ a\ GA\ bundle\ over}\ M\ {\rm along}
\ D\}$. Then ${\cal M}(M;D)$ is a unique factorization Tannakian category
and $\omega_D:{\cal M}(M;D)\to {\rm Vec}_{\bf C}$ is a fiber functor.}
\vskip 0.30cm
\noindent
{\it Proof.} (1)  By a result of Mehta-Seshadri [MS, Prop. 1.15],
${\cal M}(M;D)$ is an abelian category. Then from the unitary
representation interpretation of a GA bundle, a fundamental result due to
Seshadri, (see [MS, Thm 4.1], also in Step 3 of Section 5 below,)
${\cal M}(M;D)$ is closed under tensor product. The rigidity may be
checked directly.

\noindent
(2) Since the Jordan-H\"older graded bundle is a direct sum
of stable and hence irreducible objects and is unique, (see [MS, Rm
1.16],) so, ${\cal M}(M;D)$ is  a unique factorization category.

\noindent
(3) By definition, we know that the functor $\omega$ is exact and tensor.
So we should check whether it is faithful. This then is a direct
consequence of the fact that ${\cal M}(M;D)$ is  a unique factorization
category and that any morphism between two irreducible objects is either
zero or a constant multiple of the identity map.
\vskip 0.30cm
\noindent
{\li 3. Reciprocity Map}
\vskip 0.30cm
In $[\Sigma]$, choose its associated Jordan-H\"older graded bundle ${\rm
Gr}(\Sigma)$ as a representative. Then by the above mentioned
fundamental result of Seshadri,
${\rm Gr}(\Sigma)$ corresponds to a unitary representation $\rho_{{\rm
Gr}(\Sigma)}:\pi_1(M^0)\to U(r_\Sigma)$, where $r_\Sigma$ denotes the rank
of
$E(\Sigma):=E$.

For each element $g\in \pi_1(M^0)$, we then obtain a ${\bf
C}$-isomorphism of $E({\rm Gr}(\Sigma))|_P$. Thus, in particular, we get
a natural morphism $$W:\pi_1(M^0)\to {\rm Aut}^\otimes\omega_D.$$

Now note that $\pi_1(M^0)$ is generated by $2g$ hyperbolic
transformations $A_1,B_1,\dots,A_g,B_g$ and $N$ parabolic transformations
$S_1,\dots,S_N$ satisfying a single relation
$A_1B_1A_1^{-1}B_1^{-1}\dots A_gB_gA_g^{-1}B_g^{-1}S_1\dots S_N=1,$ and
that $\rho_{{\rm Gr}(\Sigma)}(S_i^{e_i})=1$ for all $i=1,\dots, N$.
([MS, \S 1].)  Denote by $J(D)$
the normal subgroup of $\pi_1(M^0)$ generated by
$S_1^{e_1},\dots,S_N^{e_N}$. Then naturally we obtain the following
{\it reciprocity map}
$$W(D):\pi_1(M^0)/J(D)\to {\rm Aut}^\otimes \omega_D.$$
\vskip 0.30cm
\noindent
{\li 4. Main Theorem}
\vskip 0.30cm
As usual, a Galois covering $\pi:M'\to M$ is called branched at most at
$D$ if (1) $\pi$ is  branched at $P_1,\dots,P_N$; and (2) the
ramification index $e_i'$ of points over $P_i$ divides $e_i$ for all
$i=1,\dots, N$. Clearly, by changing $D$, we get all finite Galois
coverings of $M$.
\vskip 0.30cm
\noindent
{\bf Main Theorem.} (1) (Existence and Conductor Theorem)
{\it There is a natural one-to-one correspondence $w_D$ between $$\{{\bf
S}:{\rm finitely\ completed\ Tannakian\ subcategory\ of} \ {\cal
M}(M;D)\}$$ and $$
\{\pi:M'\to M:{\rm finite\ Galois\ covering\ branched\ at\ most\ at}\
D\};$$}
\noindent
(2) (Reciprocity Law) {\it There is a natural group isomorphism}
$${\rm Aut}^\otimes (\omega_D\big|_{\bf S})\simeq {\rm Gal}\,(w_D({\bf
S})).$$
\vskip 0.20cm
\noindent
{\li 5. Proof}
\vskip 0.30cm
\noindent
{\it Step 1: Galois Theory.} By a result of Bungaard, Nielsen, Fox, M.
Kato, and Namba, (see e.g., [Na, Thms 1.2.15 and 1.3.9],) we know
that the assignment
$(\pi:M'\to M)\mapsto
\pi_*(\pi_1(M'\backslash\pi^{-1}\{P_1,\dots,P_N\}))$ gives a one-to-one
correspondence between isomorphism classes of finite Galois coverings
$\pi:M'\to M$ branched at most at $D$ and finite index (closed) normal
subgroups $K=K(\pi)$ of $\pi_1(M^0)$ containing $J(D)$. Moreover, we have
a natural isomorphism ${\rm Gal}(\pi)\simeq \pi_1(M^0)/K(\pi)\Big(\simeq
\big(\pi_1(M^0)/J(D)\big)\Big/\big(K(\pi)/J(D)\big)\Big).$
Thus the problem is transformed to the one for finite index
normal subgroups of $\pi_1(M^0)$ which contain $J(D)$, or the same, finite
index normal subgroups of $G(D):=\pi_1(M^0)/J(D)$.

\noindent
{\it Step 2: Namba Correspondence.} Consider now the category ${\bf T}(D)$
of equivalence classes of unitary representations of $G(D)$. Clearly,
${\bf T}(D)$ forms a unique factorization Tannakian category, whose fiber
functor
$\omega(D)$ may be defined to be the forget functor. Now fixed once for
all a representative
$\rho_\Sigma:G(D)\to U(r_\Sigma)$ for each equivalence classes $[\Sigma]$.
(The choice of the representative will not change the essentials below as
the resulting groups are isomorphic to each other.)

Let {\bf S} be a finitely completed Tannakian subcategory of ${\bf
T}(D)$. Then as in the definition of reciprocity map above, we have a
natural morphism $G(D)\buildrel\omega_{\bf S}\over \to {\rm Aut}^\otimes
\omega\Big|_{\bf S}.$ Denote its kernel by $K({\bf S})$. Then, by
definition, $G(D)/K({\bf S})$ is a finite group, and
$K({\bf S})=\cap_{[\Sigma]\in{\bf S}}{\rm ker}\rho_\Sigma.$

Since {\bf S} is finitely completed, there exists a finite set
$S=\{[\Sigma_1],[\Sigma_2],\dots,[\Sigma_t]\}$ which generates {\bf S}
as a completed Tannakian subcategory. Set
$[\Sigma_0]:=\oplus_{i=1}^t[\Sigma_i]$. Then  for any $[\Sigma]\in {\bf
S}$, ${\rm ker}(\rho_{\Sigma_0})\subset {\rm ker}(\rho_{[\Sigma]})$,
since {\bf S} is generated by $S$. Also, by definition,
${\rm ker}(\rho_{\Sigma_0})=\cap_{i=0}^t {\rm ker}(\rho_{\Sigma_i})$.
Thus $K({\bf S})={\rm ker}(\rho_{\Sigma_0})$.

Therefore, for any $[\Sigma]\in {\bf S}$,
$\rho_\Sigma=\tilde\rho_\Sigma\circ\Pi(D;{\bf S})$ where $\Pi(D;{\bf
S}):G(D)\to G(D)/K({\bf S})$ denotes the natural quotient map and
$\tilde \rho_\Sigma$ is a suitable unitary representation of $G(D)/K({\bf
S})$.

Now set $\tilde {\bf S}:=\{[\tilde\rho_{\Sigma}]:[\Sigma]\in {\bf S}\}$.
$\tilde {\bf S}$ is a finitely completed Tannakian subcategory in
${\cal U}ni(G(D)/K({\bf S}))$, the category of equivalence classes
of unitary representations of $G(D)/K({\bf S})$.

In particular, since the unitary
representation $\tilde\rho_{\Sigma_0}$ of $G(D)/K({\bf S})$ maps
$G(D)/K({\bf S})$ injectively into its image,  for any two
elements $g_1,g_2\in G(D)/K({\bf S})$,
$\tilde\rho_{\Sigma_0}(g_1)\not=\tilde\rho_{\Sigma_0}(g_2)$.

With this, by applying the van Kampen Completeness Theorem ([Ka]),
which claims that for any compact group $G$, if Z is a subset of the
category
${\cal U}ni(G)$ such that for any two elements $g_1,g_2$, there exists a
representation  $\rho_{g_1,g_2}$ in $Z$ such that
$\rho_{g_1,g_2}(g_1)\not=\rho_{g_1,g_2}(g_2)$, then the completed
Tannakian subcategory generated by $Z$ is the whole
category ${\cal U}ni(G)$ itself, we conclude that $\tilde {\bf S}={\cal
U}ni(G(D)/K({\bf S}))$. But as categories, {\bf S} is equivalent to
$\tilde {\bf S}$, thus by the Tannaka duality, (see e.g., [DM] or [Ta])
we obtain a natural isomorphism
${\rm Aut}^\otimes\omega|_{\bf S}\simeq G(D)/K({\bf S}).$

On the other hand, if $K$ is a finite index (closed) normal subgroup of
$G(D)$. Set $\tilde {\bf S}:={\cal U}ni(G(D)/K)$ with the fiber functor
$\omega_{\tilde {\bf S}}$. Compositing with the natural quotient map
$\Pi:G(D)\to G(D)/K$ we then obtain an equivalent  category {\bf S}
consisting of corresponding unitary representations of
$G(D)$. {\bf S} may also be viewed as a Tannakian subcategory of ${\cal
U}ni(G(D))$. We next show that indeed such an {\bf S} is a finitely
completed Tannakian subcategory.

 From definition, ${\rm Aut}^\otimes\omega|_{\bf S}\simeq {\rm
Aut}^\otimes\omega_{\tilde{\bf S}}$ which by the Tannaka duality theorem
is ismorphic to $G(D)/K$. So it then suffices to show that {\bf S} is
finitely generated. But this is then a direct consequence of the fact that
for any finite group there always exists a unitary representation such
that the group is injectively mapped into the unitary group.

\noindent
{\it Step 3: A Micro Reciprocity Law.} With Steps 1 and 2, the proof of
the Main Theorem is then completed by the following

\noindent
{\bf Weil-Narasimhan-Seshadri Correspondence} (Seshadri, [MS, Thm 4.1])
{\it There is a natural one-to-one
correspondence between  isomorphism classes of unitary
representations of fundamental groups of $M^0$
and equivalence classes of semi-stable parabolic bundles over $M^0$ of
parabolic degree zero.}

Indeed, with this theorem, the Seshadri equivalence classes
of GA bundles over
$M$ along
$D$ correspond  naturally in one-to-one to the equivalence classes of
unitary representations of the group $\pi_1(M^0)/J(D)$. Thus by Step 2,
the finitely completed Tannakian subcategories of ${\cal M}(M;D)$ are in
one-to-one correspond to the finite index closed normal subgroup of
$\pi_1(M^0)/J(D)$, which by Step 1 in one-to-one correspond to the finite
Galois coverings of $M$ branched at most at $D$. This gives the existence
theorem. Along the same line, we have the reciprocity law as well.
\vskip 0.45cm
\centerline {\li REFERENCES}
\vskip 0.25cm
\item{[DM]} P. Deligne \& J.S. Milne, Tannakian categories, in {\it Hodge 
Cycles, Motives and
Shimura Varieties}, LNM {\bf 900}, (1982), 101-228
\vskip 0.25cm
\item{[Ka]} E. van Kampen, Almost periodic functions and compact groups,
Ann. of Math. {\bf 37} (1936), 78-91
\vskip 0.25cm
\item{[MS]} V.B.  Mehta \& C.S.  Seshadri, Moduli of vector bundles on 
curves with
parabolic structures.  Math. Ann.  {\bf 248}  (1980), no. 3, 205--239.
\vskip 0.25cm
\item{[MFK]} D. Mumford, J. Fogarty \& F. Kirwan, {\it Geometric
Invariant Theory}, Springer-Verlag, 1994
\vskip 0.25cm
\item{[Na]} M.  Namba, {\it Branched coverings and algebraic functions}.
Pitman Research Notes in Mathematics Series {\bf 161}, Longman Scientific
\& Technical, 1987
\vskip 0.25cm
\item{[NS]} M.S.  Narasimhan \& C.S.  Seshadri,  Stable
and unitary vector bundles on a compact Riemann surface. Ann. of Math.
{\bf 82} (1965), 540-567
\vskip 0.25cm
\item{[Se1]} C.S. Seshadri,  Moduli of $\pi$-vector bundles over an
algebraic curve.  {\it Questions on Algebraic Varieties}, C.I.M.E.,
III, (1969) 139--260
\vskip 0.25cm
\item {[Se2]} C. S. Seshadri, {\it Fibr\'es vectoriels sur les courbes 
alg\'ebriques}, Asterisque
{\bf 96}, 1982
\vskip 0.25cm
\item{[Ta]} T. Tannaka, {\it Theory of topological groups}, Iwanami, 1949
(in Japanese)
\vskip 0.25cm
\item{[W]} A. Weil, G\'en\'eralisation des fonctions ab\'eliennes, J.
Math Pures et Appl, {\bf 17}, (1938) 47-87
\vskip 0.25cm
\item{[We]} L. Weng, A Program for Gromatric Arithmetic, preprint, 2001
\end